\documentstyle[11pt,fleqn]{article}
\begin{document}
   \renewcommand{\theequation}{\arabic {equation}}
  \newcommand{\bi}{\begin{equation}}
  \newcommand{\ei}{\end{equation}}
   \date{}
   \baselineskip 24pt
    \title { Inequalities for Continued Fractions,\ II}
\author{Zaizhao  Meng}
\maketitle \baselineskip 24pt
\baselineskip 24pt
 \begin{center}
  In this paper, we investigate the monotone property of the continued fractions
   $G(m,\lambda)$ as a function of $m$ and $\lambda$.  In particular, we obtain
     new inequality for the relative continued fractions.
\end{center}
{\bf Key Words:} continued fractions; monotone property.
$$ $$
\section*{ 1. Introduction}
 We use the symbols and the properties of continued fractions in [1].\\
Let
$$<x_{0},x_{1},x_{2},\ldots> = x_{0}+\frac{1}{x_{1}+\frac{1}{{x_{2}+}_{\ddots}}},$$
$$<x_{0},x_{1},x_{2},\ldots,x_{n}> = x_{0}+\frac{1}{x_{1}+\frac{1}{x_{2}{+}_{\ddots+\frac{1}{x_{n}}}}},$$
and $P_{-2}=0,\ P_{-1}=1,\ Q_{-2}=1,\ Q_{-1}=0,$ $
P_{n}=x_{n}P_{n-1}+P_{n-2},\\
 Q_{n}=x_{n}Q_{n-1}+Q_{n-2},
$ then, for $n\geq 0$, we have
$$<x_{0},x_{1},x_{2},\ldots,x_{n}> =\frac{P_{n}}{Q_{n}}.$$

Suppose that $m> -1,\ \lambda >0$, and let $ x_{0}=m\lambda,\ x_{1}=(m+1)\lambda,\ \ldots,\\ x_{j}=(m+j)\lambda,\ \ldots,\
x_{n}=(m+n)\lambda . $\  We define the functions
$$G_{n}(m,\lambda):=\frac{P_{n}}{Q_{n}},$$
$$G(m,\lambda):=\lim\limits_{n\rightarrow \infty}G_{n}(m,\lambda).$$
$G(m,\lambda)=<m\lambda,(m+1)\lambda,(m+2)\lambda,\ldots> = m\lambda +\frac{1}{(m+1)\lambda +\frac{1}{{(m+2)\lambda +}_{\ddots}}}.$\\
The purpose of this paper is to obtain new inequality for  $G(m,\lambda)$.\\
\ \ {\bf THEOREM .}\ \ For $m\geq 0,\ \lambda >0$ , we have
$$
G(m+1,\lambda)>\frac{m\lambda}{2}+\sqrt{\frac{m^{2}\lambda^{2}}{4}+1}>G(m,\lambda).$$
 In particular,  for $m\geq 1,\ \lambda >0$, \ \ $ G(m,\lambda)>1.$\\
 $\lambda +\frac{1}{2\lambda +\frac{1}{{3\lambda +}_{\ddots}}}>1,\ \lambda>0.$

\section*{2.\  \ Proof of Theorem }

By Corollary 2 of [1], for   $\lambda >0$, we have
$$G(m+1,\lambda)>G(m,\lambda)=m\lambda+\frac{1}{G(m+1,\lambda)},$$
consequently, we obtain $ G(m+1,\lambda)>\frac{m\lambda}{2}+\sqrt{\frac{m^{2}\lambda^{2}}{4}+1}>G(m,\lambda),$ the  theorem  follows.\\
We obtain $G(0,\lambda)=\frac{1}{G(1,\lambda)}<1$,\ for fixed
$\lambda >0$, there exists $\alpha(\lambda)\in(0,1)$ such that
$G(\alpha(\lambda),\lambda)=1$, by Theorem 3 of [1], for $m>-1, \
\lim\limits_{\lambda\rightarrow 0+}G(m,\lambda)=1$,  hence , there
exists $m\in(0,1)$, the function $G(m,\lambda)$ of $\lambda$ is not
monotonically increasing in $ (0,+\infty)$.

{\small E-mail:\ mengzzh@126.com}
 \end{document}